\documentclass[11pt]{article}
\pdfoutput=1

\usepackage{etex}
\usepackage{enumerate,amsmath,amssymb,amsthm}
\usepackage{setspace}
\usepackage{sectsty}
\usepackage{color}
\usepackage{subfigure}
\usepackage{hyperref}
\sectionfont{\large}
\subsectionfont{\normalsize}
\usepackage[compact]{titlesec}
\usepackage[margin = 1 in]{geometry}
\usepackage{graphicx}
\def\Box{\hskip 0.2cm \kern  -0.2cm\hbox{\vrule width 0.2cm height 0.2cm}}

\title{Towards a $q$-analogue of the Harer-Zagier formula via rook placements}
\author{Max Wimberley (MIT)}

\newcommand{\qbinom}[2]{\genfrac{[}{]}{0pt}{}{#1}{#2}}

\DeclareMathOperator{\inv}{inv}

\DeclareMathOperator{\sz}{size}
\DeclareMathOperator{\cn}{cn}

\begin{document}

\numberwithin{equation}{section}
\newtheorem{theor}[equation]{Theorem}
\newtheorem{prop}[equation]{Proposition}
\newtheorem{lem}[equation]{Lemma}
\newtheorem{cor}[equation]{Corollary}
\newtheorem{conj}[equation]{Conjecture}
\newtheorem*{bijTheor*}{Theorem \ref{theor:Bijection}}
\newtheorem*{bijCor*}{Corollary \ref{cor:mainresultscor}}
\newtheorem*{qConj*}{Conjecture \ref{conj:qIdent}}
\newtheorem*{qHermiteMoments*}{Theorem \ref{theor:qHermiteMoments}} 
\theoremstyle{remark}
\newtheorem{rem}[equation]{Remark}
\newtheorem{ex}[equation]{Example}
\newtheorem{exs}[equation]{Examples}
\newtheorem{questions}[equation]{Questions}

\theoremstyle{definition}
\newtheorem{definition}[equation]{Definition}

\maketitle

\begin{abstract}
In 1986 Harer and Zagier computed a certain
matrix integral to determine an influential closed-form formula for the number of
(orientable) one-face maps on $n$ vertices colored from $N$ colors. Kerov (1997)
provided a proof which computed the same matrix integral differently,
which gave an interpretation of these numbers as also counting the
number of placements of non-attacking rooks on Young diagrams. Bernardi
(2010) provided a bijective proof of this formula by putting one-face
maps in bijection with tree-rooted maps, which are orientable maps
with a designated spanning tree. In the first part of the paper, we explore the connection between
these rook placements and tree-rooted maps by developing a bijection
between these objects. Rook placements on Young diagrams have a $q$-analogue due to Garsia and Remmel (1986). 
In the second part of the paper, we propose a
statistic on rook placements that leads to a conjectured identity
which is a $q$-analogue of part of the Harer-Zagier formula. This identity is also expressed in terms
of moments of orthogonal polynomials which are rescaling of
$q$-Hermite polynomials. We then use these moments to give a
recurrence for the proposed $q$-analogue.

Note: Conjecture \ref{conj:qIdent} has been proved by Stanton \cite{DSpf}.

\end{abstract}

\section{Introduction: Overview of Problem and Main Results}

The Harer-Zagier formula involves the enumeration of \textit{unicellular (one-face) maps} -- embeddings of 
graphs with $n$ edges on orientable surfaces (up to homeomorphism) such that cutting the surface along the edges of the graph results in a 
disk (the \textit{face}). Because the map is required to have one face, the equation $2-2g = v - e + f = v - n + 1$ shows that the genus $g$ is 
determined uniquely by the number $v$ of vertices and vice-versa. The formula is
\begin{equation}\label{eq:HZForm}
 C(n,N) := \sum_{g=0}^{\lfloor n/2 \rfloor}\varepsilon_g(n)N^v = (2n-1)!!\sum_{k\geq0}\binom{N}{k+1}\binom{n}{k}2^{k},
\end{equation}
where $\varepsilon_g(n)$ is the number of distinct unicellular maps with $n$ edges on a surface of genus $g$ (for definitions, see
Section \ref{sec:bkgrnd}) and $N$ is a positive integer. Equation \eqref{eq:HZForm} was demonstrated by computing the matrix integral
\begin{equation}\label{eq:MatrixInt}
 C(n,N) = \int_{H_N}\!\textrm{tr}\, Z^{2n}\, \mu_N(\textrm{d}Z)
\end{equation}
over the space $H_N$ of $N\times N$ Hermitian matrices $Z = (z_{jk})$ with $z_{jk} = x_{jk} + iy_{jk}$ with the appropriate Gaussian measure 
$\mu_N(dZ)$ \cite[Eq. (1.2)]{K}.

\begin{figure}[h]
\centering
\includegraphics[scale=0.9]{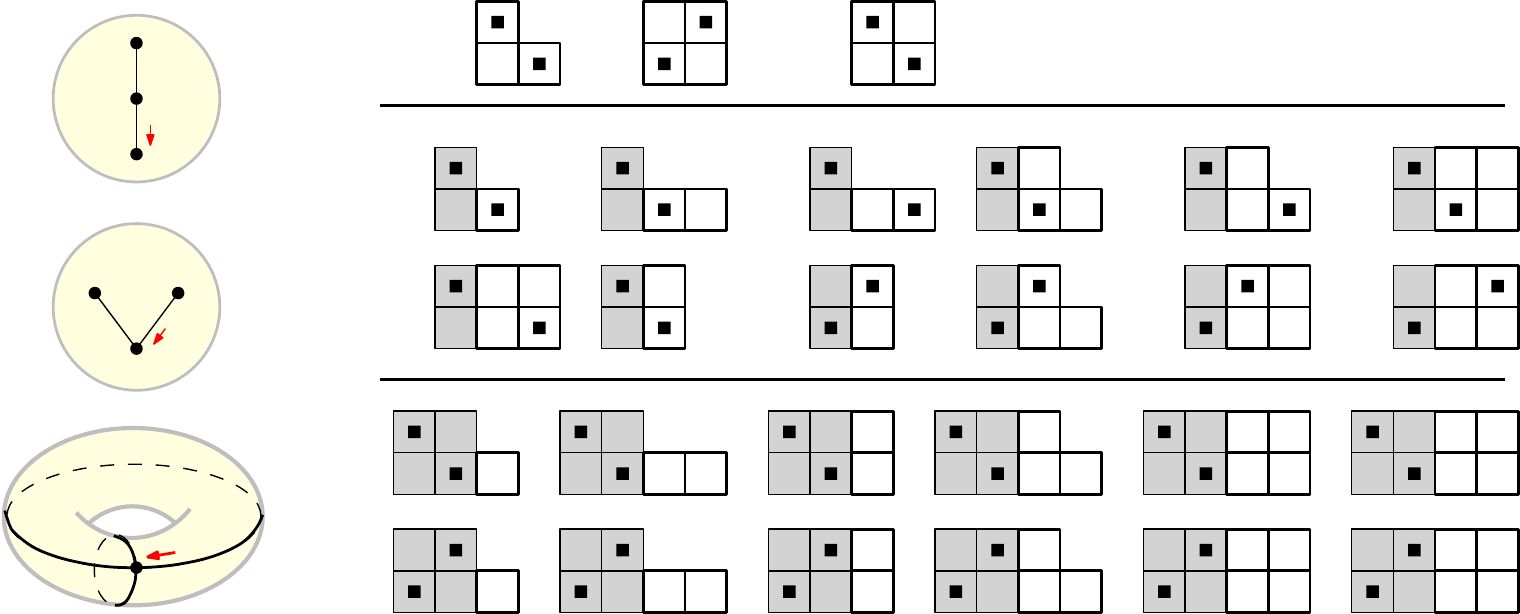}
\[
2N^3 + N = 3\binom{N}{1} + 12\binom{N}{2} + 12\binom{N}{3}.
\]
\caption{Example of the Harer-Zagier formula \eqref{eq:HZForm} for $n=2$ and the Kerov's
correspondence with rook placements.  Of the three
  rooted one-face maps with $n=2$ edges on the left, two of them are embedded on the sphere ($g=0$) and one is
  embedded on the torus ($g=1$). For the rook placements on
  the right, there
are $3, 12$ and $12$ placements of $n=2$ rooks on Young diagrams inside a $2
\times 2$ board with $k=0,1$ and $2$ columns of length $2$ added to
the left that are required to have a rook.}
\label{fig:exKerov}
\end{figure}

 Kerov \cite{K} showed that the matrix
integral (\ref{eq:MatrixInt}) could also be interpreted as counting ways of placing $n$ non-attacking rooks on a board consisting
of any Young diagram that can fit in an $n\times n$ board with $k$
columns of length $n$ added on the left that are required each to have
a rook (see Figure~\ref{fig:exKerov}). To do this,
he first transformed the integral into a sum of moments of Hermite polynomials, then evaluated them via the Flajolet-Viennot
combinatorial theory of orthogonal polynomials \cite{GV}, and also put the rook placements in bijection with certain involutions
to enumerate them explicitly.

Combinatorial proofs of \eqref{eq:HZForm} have been given by
Lass \cite{
L}, Goulden-Nica \cite{GN}, and Bernardi \cite{OB}. In particular, Bernardi proved \eqref{eq:HZForm} by interpreting each side as counting all colored unicellular
maps with colors chosen from $[N] = \{1,2,\ldots,N\}$ (we will use the notation $[n] = \{1,2,\ldots,n\}$ throughout), where on the left side 
each vertex of
the map is colored independently with a color from $[N]$ and on
the right side a subset of $k+1$ colors is first chosen and the map is required to use all of them. The right side was counted directly
via a bijection between unicellular maps with vertices colored with colors in $[k+1]$ and \textit{tree-rooted} maps with vertex set $[k+1]$,
maps whose graphs contain distinguished \textit{spanning trees}, each of which has a \textit{root vertex}.

A natural question arises from these results. Mainly, the only extant connection between unicellular maps or
(equivalently) Bernardi's tree-rooted maps and Kerov's rook placements is through the matrix integral. To provide a more direct connection between
these two types of objects, we establish a bijection between tree-rooted maps and involutions which, along with Kerov's and Bernardi's
bijections, yields the following corollary:
\begin{bijCor*}
 There is an explicit bijection mapping rook placements in $RC_k(n,s)$ to tree-rooted maps with $n$ edges 
 and vertex set $\{s+1\}\cup V$ where $V\subset [s]$ and $\#V = k$.
\end{bijCor*}
The set $RC_k(n, s)$ of rook placements is defined
below. We note that Bernardi~\cite{OBpc} also has an explicit bijection between the
objects above. 

One motivation for finding such a bijection is that it may preserve
interesting properties of either object. In particular, from the work
of Garsia-Remmel \cite{GR} rook
placements on Young diagrams have very simple $q$-analogues. One might then ask what could serve as a $q$-analogue of the numbers $C(n,N)$, and whether this $q$-analogue can be expressed in other forms, in the same
way that $C(n,N)$ can be expressed as a matrix integral or as a sum of moments of Hermite polynomials. Section~\ref{sec:qAnalogues} addresses this with the following conjecture (verified
up to $n=10$ and $s=5$):
\begin{qConj*}
The following identity holds:
\begin{align}
\sum_{\mu \in Y(n,s)} \prod_{i=1}^{n}q^{\mu_i-i}[\mu_i-i+1]_q &= \sum_{\mu \in Y(n,s)}\sum_{\substack{\textrm{rook placements}\\ C \,\textrm{on } \mu}}q^{\textrm{inv}(C) + |\mu| - \binom{n+1}{2}} \\
&=\sum_{k\geq 0} q^{n(s-k) + \binom{k}{2}}[2n-1]_q!!\qbinom{s}{k}_q\qbinom{n}{k}_q\prod_{i=1}^k(1 + q^{n+i}),
\end{align}
\end{qConj*}
where $Y(n,s)$ is the set of Young diagrams $\mu$ with $n$ rows $\mu_i$ such that $s \leq \mu_1 \leq \mu_2 \leq \ldots \leq \mu_n \leq n+s$.
In Section~\ref{sec:orthoPolys}, we show that the expressions in Conjecture~\ref{conj:qIdent} are indeed
related to a $q$-analogue $H_n(x,q)$ of the Hermite Polynomials defined by $H_{n+1}(x,q) = xH_n(x,q) - q^{n-1}[n]_qH_{n-1}(x,q)$
with $H_0 = 1$ and $H_1 = x$. This relationship is described in the following theorem:
\begin{qHermiteMoments*}
The moments of the $q$-Hermite polynomials $H_n(x,q)$ against the Gaussian $q$-distribution $w_G(x) = E_{q^2}^{\frac{-q^2x^2}{[2]_q}}$ are given by
\[
\frac{1}{q^{\binom{s}{2}}[s]_q!c(q)}\int_{-\nu}^{\nu}x^{2n}H_s^2(x,q)w_G(x)d_qx = \sum_{\mu \in Y(n,s)}\prod_{i=1}^nq^{\mu_i - i}[\mu_i - i + 1]_q,
\]
where
\[
 \nu = \frac{1}{\sqrt{1-q}} \quad \textrm{and} \quad c(q) = 2(1-q)^{\frac{1}{2}}\sum_{m=0}^{\infty}\frac{(-1)^mq^{m(m+1)}}{(1-q^{2m+1})\prod_{i=0}^{m-1}(1-q^{2i+2})}.
\]
\end{qHermiteMoments*}
(For definitions of the $q$-number $[n]_q$, $q$-integration and the Gaussian $q$-distribution, see Section~\ref{sec:orthoPolys}).

Work towards a proof of Conjecture~\ref{conj:qIdent} is given in Sections~\ref{sec:orthoPolys}~and~\ref{sec:qAnalogues}, in the form of
a detailed study of these $q$-Hermite polynomials and identities involving them and their moments (of the form in Theorem~\ref{theor:qHermiteMoments}). We note that Stanton has 
provided a proof of Conjecture~\ref{conj:qIdent} \cite{DSpf} using techniques involving hypergeometric series which can be found in \cite{GRhs}, but that finding a combinatorial proof is still 
an open problem.

\noindent \textbf{Outline} \quad In short, the structure of the paper is as follows: In Section \ref{sec:bkgrnd}, we give background 
material and definitions we give as much background about maps and rook placements as is needed to discuss the bijection in 
Section \ref{sec:Bijection}. In Section \ref{sec:Bijection}, we describe and demonstrate the veracity of our bijection. In Section \ref{sec:orthoPolys},
we discuss the form of $q$-Hermite polynomials that we use in this paper (as well as supplementing necessary background material
about $q$-analogues). In Section \ref{sec:qAnalogues}, we
present a conjectured identity that gives a $q$-analogue for a formula counting rook placements, and use the integral 
formulation to provide a recurrence for this quantity. We conclude with closing remarks and possible future work in Section \ref{sec:conclusion}.

\noindent \textbf{Acknowledgements} \quad 
I would like to thank everyone who made this project (which was part of the MIT SPUR program in the Summer of 2012) possible -- Pavel Etingof, the head of SPUR;
Olivier Bernardi and my mentor Alejandro Morales for suggesting and designing the project; Alexei Borodin for discussing the problems and sponsoring the research over the Fall of 2012; 
and everyone who discussed these problems with us: Guillaume Chapuy, Alan Edelman, Praveen Venkataramana, Jang Soo Kim, Dennis Stanton, Jonathan Novak, and Richard Stanley.

\section{Background and Definitions}\label{sec:bkgrnd}
We begin by reviewing maps, which constitute one of the types of objects in the bijection in Section~\ref{sec:Bijection}. Maps are treated extensively in \cite{LZ}.

A \textit{unicellular map} is an embedding of a connected graph in a smooth, compact surface (up to homeomorphism) with the property that the complement of the graph in the surface
is homeomorphic to a disk. A \textit{rooted} unicellular map is a unicellular map with one side of an edge distinguished as the root and oriented. For our purposes, the surface containing a map
must be orientable. Rooted unicellular maps with $n$ edges can be obtained by gluing the edges of a rooted $2n$-gon ($2n$-gon with one edge distinguished as the root and assigned a direction). 
A gluing of the $2n$-gon has $n$ edges and 1 face by construction, so the formula for the Euler characteristic $\chi$ of the surface in which the
map is embedded is $\chi = v - n + 1$, where $v$ is the number of vertices. Since $\chi = 2 - 2g$ where $g$ is the genus of the surface, we have $2 - 2g = v - n + 1 \rightarrow
v = n + 1 - 2g$. So we can take the sum $\sum_{g = 0}^{\lfloor n/2 \rfloor} \varepsilon_g(n)N^v = \sum_{g = 0}^{\lfloor n/2 \rfloor} \varepsilon_g(n)N^{n + 1 - 2g}$ where $\varepsilon_g(n)$ 
is the number of distinct unicellular maps with $n$ edges on a surface of genus $g$, and $N$ is a positive integer. This sum
can be interpreted as counting the unicellular maps with $n$ edges and vertices colored from $[N]$. 

Because multiple vertices may have the same color, we can count the same objects by first choosing a set of $k+1$ colors (in $\binom{N}{k+1}$ ways) and insisting
the coloring use all of them. In \cite{OB}, Bernardi bijectively counted the number of unicellular maps with $n$ edges that use exactly $k+1$ colors by constructing
a bijection to \textit{tree-rooted} maps, which are maps with a designated rooted spanning tree $T$ (a spanning tree of a map is a subgraph of the underlying graph 
of the map that includes all of the vertices and is a tree; to say that it is rooted means that it has a vertex distinguished 
as the root). By construction, if a unicellular map has $n$ edges and uses $k+1$ colors, the corresponding tree-rooted map
has $n$ edges and vertex set $[k+1]$.

\begin{figure}[h]
\centering
\includegraphics[scale=0.9]{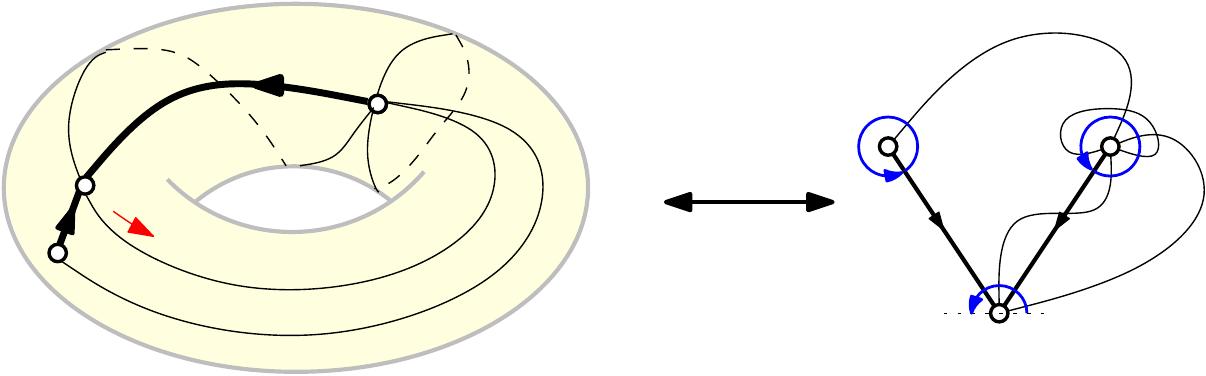}
\caption{On the left is a tree-rooted map with 3 vertices and 6 edges. On the right is the corresponding graph with its rotation system indicated by the gray arrows. The order
around the root is broken at the root-edge. The spanning trees are in bold. Bold black arrows point to the root.}
\label{fig:graphRotSysEx}
\end{figure}

An important part of Bernardi's bijection is the fact that a graph and its \textit{rotation system} (the cyclic order
of half-edges around the vertices of the graph) are sufficient to
construct a unique map (see \cite{LZ}). From this it is not difficult to see that a graph with a rooted rotation system (one vertex is the root, around which there is a total
order of half-edges) can be used to construct a rooted map (the root outgoing edge is the smallest element around the root vertex of the graph). For 
an example of a graph with a rooted rotation system and rooted spanning tree, see figure~\ref{fig:graphRotSysEx}.

We now turn to the other type of objects in the bijection in Section~\ref{sec:Bijection}, rook placements on Young diagrams.

A \textit{board} is a set of pairs $(i,j)$ which give the coordinates of the positions on the board. A non-attacking \textit{rook placement} $C$ on a board $B$ is
a set of positions $C = \{(x_i,y_i)\}\subset B$ with at most one rook in each row and column, that is, $x_i = x_j$ if and only if $i=j$ (similarly $y_i = y_j$ if and only if $i = j$) for all $(x_i,y_i)$ and $(x_j, y_j)$
in $C$. A \textit{Young diagram} with $n$ rows is a sequence of nonnegative integers $\{\mu_i\}_{i=1}^n$ such that $0 \leq 
\mu_1 \leq \mu_2 \leq \ldots \leq \mu_n$, where $\mu_i$ gives
the length of row $i$ (that is, a Young diagram can be interpreted as a board with rows of lengths given by the $\mu_i$ in order
from the top). The number of placements of $n$ rooks on a Young diagram, call it $\mu$, with $n$ rows is
\begin{equation}\label{eq:rpEnum}
 \#\{n\textrm{-rook placements } C \textrm{ on } \mu\} = \prod_{i=1}^n (\mu_i - i + 1),
\end{equation}
which follows from the fact that there are $\mu_1$ ways to place a rook on the first row, $\mu_2 - 1$ ways to place one on the second row (since $\mu_1 \leq \mu_2$),
$\mu_3 - 2$ ways to place a rook on the third row (since $\mu_1 \leq \mu_2 \leq \mu_3$), etc \cite[Section 2.3]{EC1}.

\begin{figure}[h]
\centering
\includegraphics[scale=0.9]{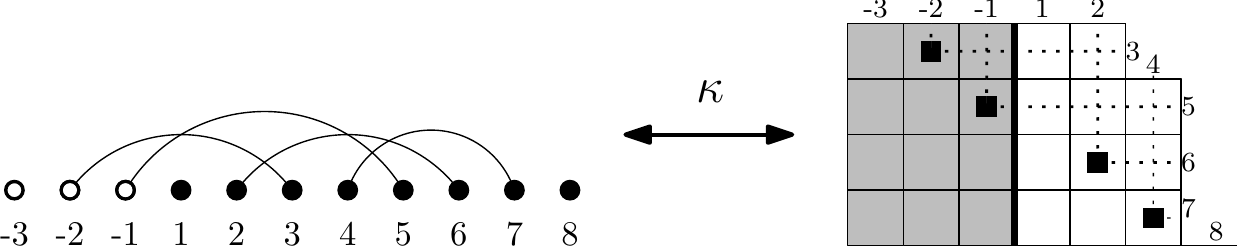}
\caption{On the left is the involution $(-3)(-2,3)(-1,5)(1)(2,6)(4)(7,8)$. In this case, $n=4$, $s=3$, and $k=2$. On the right is the corresponding rook placement. 
The numbers indicate which segments of the upper and right boundaries correspond to which points of the involution, the dotted lines show how the arcs of the involution correspond to rooks.}
\label{fig:RCIKex}
\end{figure}

We define the set $Y(n,s) = \{\textrm{boards } \mu \textrm{ with } n \textrm{ rows} \mid s \leq \mu_1 \leq \mu_2 \leq \ldots \leq \mu_n \leq s+n\}$. Each board
$\mu$ in $Y(n,s)$ can be thought of as being obtained from a Young diagram that fits in an $n\times n$ box with $s$ extra 
columns added to the left (for an example, see
Figure \ref{fig:RCIKex}. We then define the set of all rook placements of $n$ rooks on boards in $Y(n,s)$,
\[RC(n,s) = \{(\mu, C) \mid \mu \in Y(n,s), C \textrm{ a placement of } n \textrm{ non-attacking rooks on } \mu\}.\]
The subset $RC_k(n,s)\subset RC(n,s)$ consists of all elements $(\mu,C)$ of $RC(n,s)$ where exactly $k$ rooks of $C$ are in the 
first $s$ columns of $\mu$. Kerov provided
a bijection (which we will call here $\kappa$) between $RC_k(n,s)$ and another set $I_k(n,s)$, defined as the set of 
involutions (matchings) on the set $\{-s,-s+1,\ldots,-1,1,2,\ldots,2n\}$
such that exactly $k$ negative points are matched ($s-k$ are not), each negative point is matched to a positive one, and exactly 
$k$ positive points are unmatched.
Under the bijection $\kappa$ between $RC_k(n,s)$ and $I_k(n,s)$, each rook in a placement $C$ on a board $\mu$ defines an arc of the corresponding involution. One can think
of each point in the set $\{-s,-s+1,\ldots,-1,1,2,\ldots,2n\}$ as corresponding to a segment on the path defining the upper and right boundaries of $\mu$, which
consists of $s+2n$ segments ($s$ horizontal segments for the first $s$ columns, $2n$ segments defining the $n\times n$ Young diagram), so that a rook defines an arc
between the points corresponding to the horizontal segment at the top of the rook's column and the vertical segment at the right boundary of the rook's row. For an
example, see Figure~\ref{fig:RCIKex}.


\section{A Bijection Between Involutions and Tree-Rooted Maps}\label{sec:Bijection}
We define the set $I_k(n)$ as the set of involutions on the set of points $\{-k,-k+1,\ldots,-1\}\cup[2n]$ with the property that for any
involution $A \in I_k(n)$, if $p<0$ then $A(p)>0$, and $\#\{p\in[2n]\, | \, A(p) = p\} = k$ (this is essentially equivalent to
$I_k(n, s)$ where $s = k$). We now state the main theorem of this section:
\begin{theor}\label{theor:Bijection}
 There is an explicit bijection $\Psi$ mapping involutions in $I_k(n)$ to tree-rooted maps in $\mathcal{T}_n(k+1)$.
\end{theor}
Which gives the following corollary:
\begin{cor}\label{cor:mainresultscor}
 There is a bijection $\Psi'$ mapping rook placements in $RC_k(n,s)$ to tree-rooted maps with $n$ edges and vertex set $\{s+1\}\cup V$ where $V\subset [s]$ and $\#V = k$.
\end{cor}
We will prove that Corollary~\ref{cor:mainresultscor} follows from Theorem~\ref{theor:Bijection} here, while the proof of
Theorem~\ref{theor:Bijection} will be given later.

\noindent\textit{Proof of Corollary~\ref{cor:mainresultscor}.} Given a rook placement $R$ in $RC_k(n,s)$, we construct the rook placement $R'$ in $RC_k(n,k)$ by simply removing the columns
in the first $s$ columns which are empty (those which do not contain a rook). We can then find an involution $A = \kappa(R')$ 
in $I_k(n)$ corresponding bijectively to $R'$. Furthermore, define the set $V'$ as the set of
indices of columns in the first $s$ columns of $R$ which contain rooks; that is, if the $1^{\textrm{st}}$ column from the left has a rook, 1 is in the set $V'$, if the $2^{\textrm{nd}}$
column from the left has a rook, 2 is in $V'$, etc. We find $\Psi(A)$, which has vertex set $[k+1]$. There is a unique order-isomorphic correspondence between
$V'\cup\{s+1\}$ and $[k+1]$ which is obtained by writing $V' = \{i_1, i_2, \ldots, i_k\}$ with $0 < i_1 < \cdots < i_k < s + 1$ (so that $i_j \in V'$ is mapped to $j \in [k + 1]$ and $s + 1$ is mapped to 
$k + 1$); we define $\Psi'(R)$ as the relabeling of $\Psi(A)$ under this correspondence. To find $\Psi'^{-1}(M)$ for a
tree-rooted map $M$ with $n$ edges, and $k$ vertices labeled from $[s]$ as well as one vertex labeled $s+1$ (that is,
we set $s+1 = \max \{\textrm{labels of vertices}\}$) first create the map $M'$ by relabeling the vertices of $M$ with the elements of $[k+1]$ in an order-isomorphic way,
then find $\Psi^{-1}(M')$, an involution in $I_k(n)$. Next, find the rook placement $R' = \kappa^{-1}(\Psi^{-1}(M'))$.
Lastly, create $R = \Psi'^{-1}(M)$ by inserting empty columns between the first $k$ columns of $R'$ so that the indices of the columns containing rooks (which are in the first $s$ columns) match the 
original labels of $M$. \qed

Before we can prove Theorem~\ref{theor:Bijection}, we note some classical results. In our description of the bijection, we use the term \textit{path} (not the same as the
weighted paths used to compute the moments of the Hermite polynomials) to refer to a sequence of
\textit{steps}, each of which can be \textit{up} or \textit{down}. An up step begins at a level $l$ and ends at level $l+1$, whereas a down step begins
at a level $l$ and ends at level $l-1$. We can also think of these as paths from $(0, 0)$ to $(L, 0)$ (for some $L \in \mathbb{Z}^+$) in the lattice $\mathbb{Z}^2$, where the level corresponds to the 
$y$-coordinate, in which case an up step moves in the direction $(1, 1)$ while a down step moves in the direction $(1, -1)$. The first step of a path always begins at level 0 and, for our purposes, the 
last step must end at level 0 (note that this means there must be an equal number of up and down steps, hence the length of the path must be even). The number of steps in the path is the \textit{length} of 
the path. A \textit{Dyck path} is a path with no step ending on a level below 0. We will also define the number of \textit{flaws} in a path to be
the number of up steps starting below 0.

\begin{figure}[h]
\centering
\includegraphics[scale=0.9]{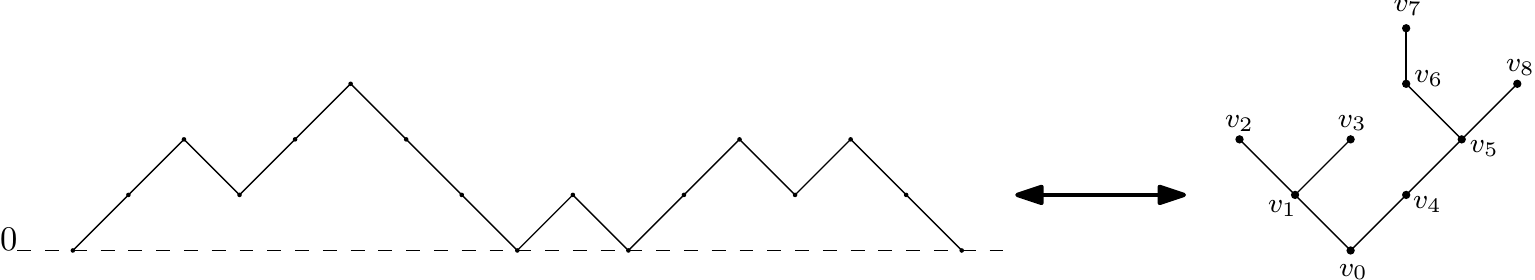}
\caption{$P$ as in Example \ref{ex:pathTree}, and the resulting tree. Note that the subscripts of the $v_i$ directly indicate
the correspondence described in Remark~\ref{rem:correspondence}, with $v_i$ corresponding to the $i^\textrm{th}$ down step.}
\label{fig:pathTreeEx}
\end{figure}

We make use of a classical bijection (which we'll denote by $\beta$) between Dyck paths and trees that uses a \textit{breadth-first} (as opposed to the also-common 
depth-first) method. Under this correspondence, the number of consecutive
up steps after each down step is the number of children of the next vertex in the tree, where the vertices are ordered lexicographically
such that the subtree beginning from any vertex $v$ comes before the next vertex to the right of $v$ on the same level. We state this
as a lemma, noting that the result is classical:
\begin{lem}\label{lem:pathTree}
 The map $\beta$ described above is bijective.
\end{lem}
\begin{rem}\label{rem:correspondence}
 The above method of constructing a tree from a Dyck path gives a natural correspondence between the down steps of the path and the
non-root vertices of the tree, where a down step $s_i$ corresponds to the vertex with (the number of) children given by the number of
up steps following $s_i$.
\end{rem}
\begin{ex}\label{ex:pathTree}
 Let $P = UUDUUDDDUDUUDUDD$. Proceeding as described above, we see that $v_0$ has two children because there are two consecutive
 up steps at the beginning. The next two down steps both have 0 up steps after them so the two children of $v_1$ have no children.
 After the next down step, there is an up step, so the next vertex $v_4$ (the right child of $v_0$) has one child, and so on. This
 path and the resulting tree are shown in Figure \ref{fig:pathTreeEx}.
\end{ex}

It is well-known that the number of Dyck paths of length $2k$ is the $k^{\rm th}$ Catalan number ${\rm Cat}_k = \frac{1}{k+1}\binom{2k}{k}$, while
the number of paths of length $2k$ is simply $\binom{2k}{k}$. This suggests a $(k+1)$-to-1 correspondence between paths and
Dyck paths. The classical Chung-Feller Theorem \cite{JR} gives the equidistribution of the number of flaws in a path; the following lemma uses
this to construct such a correspondence.
\begin{lem}\label{lem:ChungFeller}
 (Chung-Feller). The set of paths of length $2k$ can be partitioned into disjoint subsets of size $k+1$ such that for each $j$,
$0 \leq j \leq k$ there is exactly one path with exactly $j$ flaws in each subset. So paths of length $2k$ can be put into a $(k+1)$-to-1
correspondence with Dyck paths, by letting each path correspond to the Dyck path in the same subset of this partition.
\end{lem}
\noindent\textit{Proof.} We will define a set of bijections $\phi_j(P)$ between paths with $j-1$ flaws and $j$ flaws ($j > 0$). Let $P$ be a path of
length $2k$ with $j > 0$ flaws. Call the $i^{\textrm{th}}$ step of $P$ $s_i$; we
write $P = s_1s_2\ldots s_{2k}$. To find the image of $P$ under $\phi_j$ (which we want to be a path with $j-1$ flaws), find the first up step $s_f$ ending
at level 0. $\phi_j(P)$ is then given by $s_{f+1}s_{f+2}\ldots s_{2k}s_f s_1 s_2 \ldots s_{f-1}$. That this is a bijection with the desired properties is shown in
\cite{JR}. We then construct a set for each Dyck path $D$. Define $S_i(D) = \phi_i^{-1}(S_{i-1}(D))$ with $S_0(D) = D$. We then form the set $S(D) = \{S_0(D), S_1(D), \ldots, S_k(D)\}$.
which clearly contains $k+1$ elements, with $S_i(D)$ having $i$ flaws. Furthermore, because the $\phi_i$ are bijections, $S_i(D_1) = S_i(D_2)$ if and only if
$D_1 = D_2$, so that $S(D_1) \cap S(D_2) = \emptyset$ if $D_1 \neq D_2$. Since there are $\textrm{Cat}_k = \frac{1}{k+1}\binom{2k}{k}$ such
paths $D$,
\[
\#\bigsqcup_{\substack{\textrm{Dyck paths } \\ D \textrm{ of length } 2k}} S(D) =(k+1)\textrm{Cat}_k = \binom{2k}{k} = \#\{\textrm{Paths of length } 2k\},
\]
so $\bigsqcup_{D} S(D) = \{\textrm{Paths of length } 2k\}$. \qed

We are now ready to prove Theorem~\ref{theor:Bijection}.

\noindent\textit{Proof of Theorem~\ref{theor:Bijection}.} We construct the bijection $\Psi:I_k(n) \rightarrow \mathcal{T}_n(k+1)$ as follows. Let $A$ be an involution in $I_k(n)$. We indicate by $a_i$ the $i^{\textrm{th}}$
point of $A$ that is either a fixed point of $A$ or the image of a negative point under $A$. We construct a path $P = s_1s_2\ldots s_{2k}$
with labeled down steps (separate from the notation $s_i$, which simply refers to the $i^{\rm th}$ step)
from this sequence by the rule that $s_i$ is an up step if $A(a_i) = a_i$, otherwise it is a
down step with the label $p$ where $A(p) = a_i$ ($p<0$). This path starts and ends on the same level (which we'll call 0) by
virtue of the fact that there are $k$ fixed points and $k$ negative points so the numbers of up and down steps in the path are equal.

\begin{figure}[h]
\centering
\includegraphics[scale=0.9]{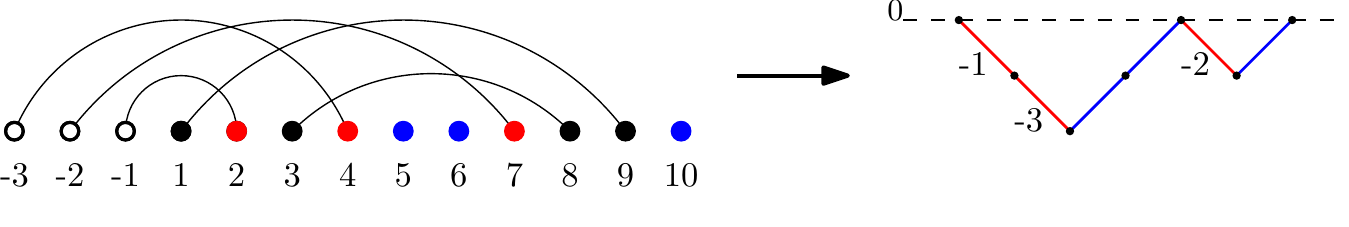}
\caption{A graphical depiction of the objects in Example \ref{ex:bijectionExample}. The coloring indicates the correspondence between
the points in $A$ and steps in $P$, and the down steps are labeled by the negative points paired with the corresponding positive
points.}
\label{fig:bijectionExamplePart1}
\end{figure}

\begin{figure}[h]
 \centering
 \includegraphics[scale = 0.75, trim = 0in 0.6in 0in 0.4in]{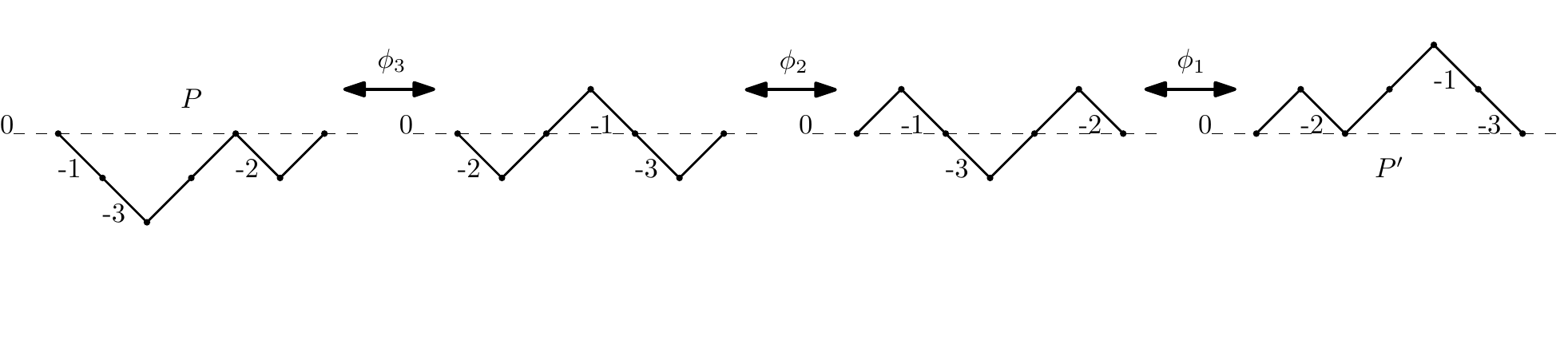}
 \caption{Application of lemma \ref{lem:ChungFeller}. The path on the far left has 3 flaws, so applying $\phi_3$ yields the
 next-leftmost path with 2 flaws, and so on, until we have a Dyck path.}
 \label{fig:chungfellerex}
\end{figure}

The number of flaws in $P$ determines
the root vertex in $\Psi(A)$: if there are $j$ flaws, the root is the vertex $j+1$.
\begin{ex}\label{ex:bijectionExample}
 Let $A \in I_3(5)$ be given by $A = (-3, 4)(-2, 7)(-1, 2)(1, 9)(3, 8)(5)(6)(10)$ (see Figure~\ref{fig:bijectionExamplePart1}). Then $a_1,\ldots,a_6 = 2, 4, 5, 6, 7, 10$ and $P = D_{-1}D_{-3}UUD_{-2}U$ where
$D_p$ indicates a down step with label $p$ and $U$ indicates an up step. This path has 3 flaws, so the root of $\Psi(A)$ will be vertex 4.
\end{ex}
We now find the Dyck path $P' = s_1's_2'\ldots s_{2k}'$ corresponding to $P$ under the correspondence described in Lemma~\ref{lem:ChungFeller}. From
this we construct the spanning tree of $T$ of $\Psi(A)$ by the procedure described in Lemma~\ref{lem:pathTree}. 
\begin{ex}\label{ex:bijectionExamplePart2}
 With $P$ and $A$ as in Example \ref{ex:bijectionExample}: $P$ has 3 flaws, so the root of the corresponding tree is vertex 4. 
 To find the tree, applying Lemma \ref{lem:ChungFeller}, we take $P' = \phi_1\phi_2\phi_3(P)$. $P'$ is a Dyck path, so we can 
 then construct a tree using $\beta$ as in Lemma \ref{lem:pathTree}, taking $T = \beta(P')$. Figure \ref{fig:chungfellerex} shows 
 the application of Lemma \ref{lem:ChungFeller} to get a Dyck path.
\end{ex}

At this point, the root vertex and the spanning tree of $\Psi(A)$ have been determined; what remains is to determine the labeling of the non-root vertices, to
insert half-edges which are not part of the spanning tree, and to pair these half-edges.

\begin{figure}[h]
\centering
\includegraphics[scale=0.9, trim= 0in 0.9in 0in 1.3in]{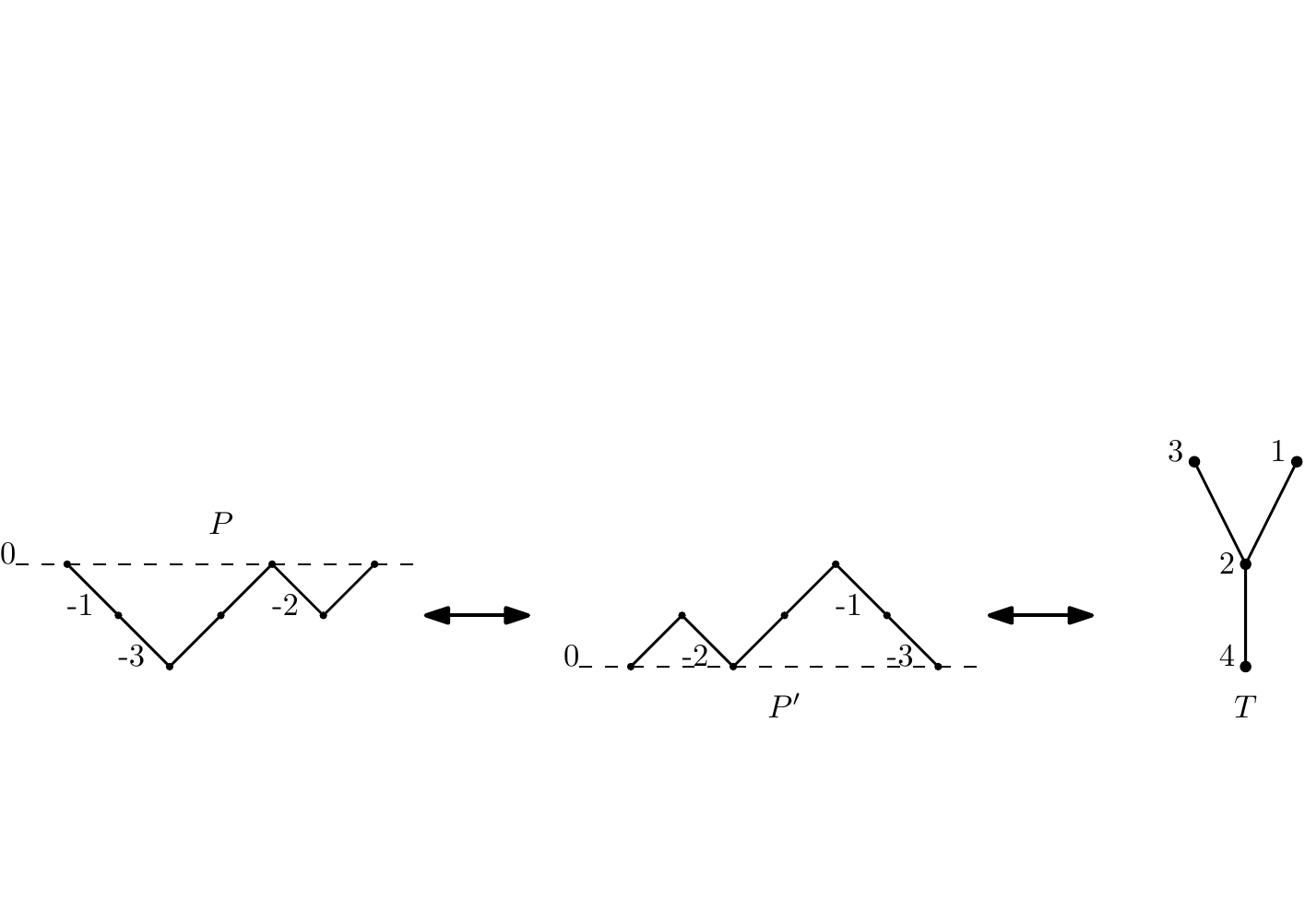}
\caption{The order-isomorphic correspondence between $\{-3,-2,-1\}$ and $\{1,2,3\}$ is $-3\leftrightarrow 1$, $-2\leftrightarrow 2$, $-1\leftrightarrow 3$. The child of 4 in the tree
corresponds to the down step labeled -2 so that vertex is labeled 2. Its left child corresponds to the down step labeled -1 so it is labeled 3, and its right child corresponds to the 
down step labeled -3 so it is labeled 1.}
\label{fig:bijectionExamplePart2}
\end{figure}

We label the non-root vertices
of $T$ as follows: As noted in Remark~\ref{rem:correspondence}, the method for building the spanning tree from $P'$ gives a natural correspondence between the down steps of $P'$ and the non-root vertices
of $T$. There is also a unique order-isomorphic (using the order on the negative points of $A$, which are the labels of the down steps)
relabeling of the down steps with the labels that must be given to the non-root vertices in $T$ (the set $[k+1]\backslash\{j+1\}$, where $j+1$ was the label given to the root).
The labeling of a non-root vertex in $T$ is then given by the label of
the corresponding southeast segment of $P'$.
\begin{ex}\label{ex:bijectionExamplePart3}
 With $A$, $P$, $P'$, and $T$ as in Examples \ref{ex:bijectionExample} and \ref{ex:bijectionExamplePart2}, we use the order-isomorphic
 correspondence $-3 \leftrightarrow 1$, $-2 \leftrightarrow 2$, and $-1 \leftrightarrow 3$ and the correspondence between the
 down steps of $P'$ and the non-root vertices of $T$ as noted in remark \ref{rem:correspondence} to label the non-root vertices.
 This is shown in Figure \ref{fig:bijectionExamplePart2}.
\end{ex}

\begin{figure}[h]
\centering
\includegraphics[scale=0.9, trim = 0in 0in 0in 0.2in]{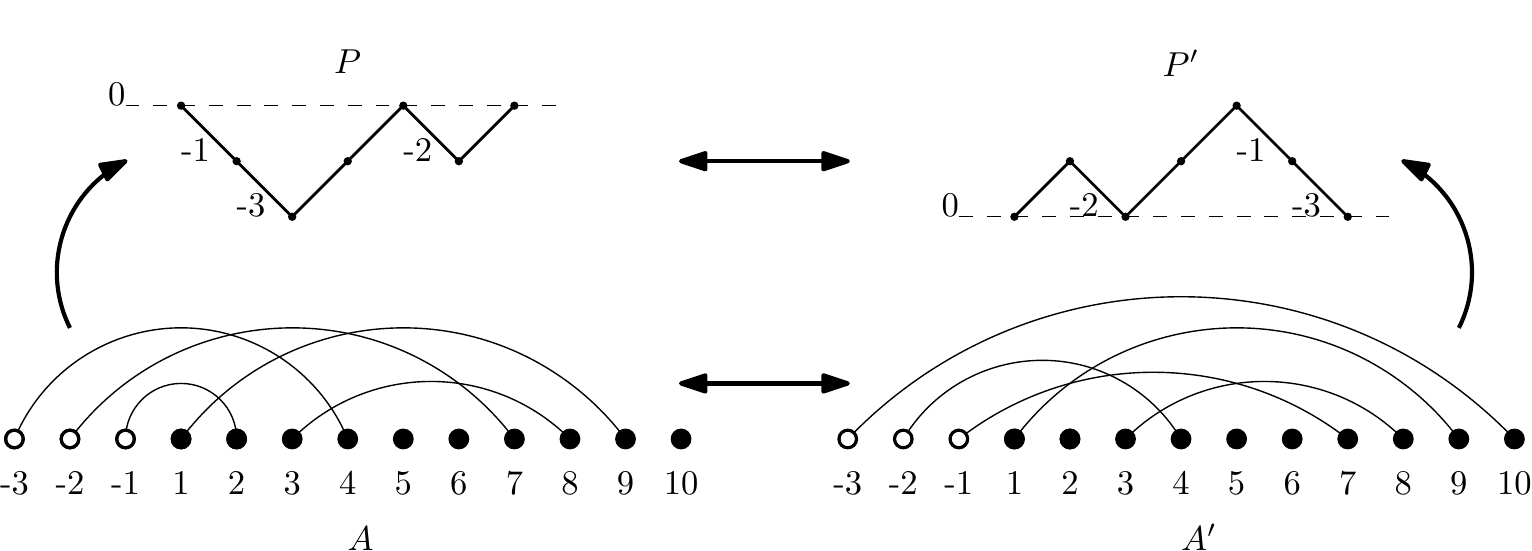}
\caption{$\sigma$ is constructed so that taking $\sigma A\sigma^{-1}$ performs the exact same rearrangement of the points in $A$ that 
we performed when rearranging the steps in $P$.}
\label{fig:bijectionExamplePart3}
\end{figure}

 To insert and pair half-edges in $T$, we first construct a modified version of $A$ from 
 the path $P'$ which simply rearranges the $a_i$ in an analogous manner:
define $A' = \sigma A\sigma^{-1}$ where $\sigma$ is the unique permutation defined
by $\sigma(a_i) = a_j$ if and only if $s_i = s_j'$. 
\begin{ex}\label{ex:bijectionExamplePart3}
 Continuing with $A$, $P$, and $P'$ as in Examples \ref{ex:bijectionExample}, \ref{ex:bijectionExamplePart2}, and \ref{ex:bijectionExamplePart3},
 we have $s_1 = s_5', s_2 = s_6', s_3 = s_3', s_4 = s_4', s_5 = s_2'$, and $s_6 = s_1'$, so $\sigma = (a_1,a_5,a_2,a_6)(a_3)(a_4) = (2,7,4,10)(5)(6)$, 
 giving $A' = (-3,10)(-2,4)(-1,7)(1,9)(2)(3,8)(5)(6)$, as desired. This is shown in Figure \ref{fig:bijectionExamplePart3}.
\end{ex}

We take the set $H$ of positive points $H=[2n]\backslash\{a_i\}_{i=1}^{2k}$, which can be partitioned into
blocks $B_i = \{p\in[2n] \, | \, a_{i-1} < p < a_i\}$ with $B_1 = \{p\in[2n] \, | \, p < a_1 \}$ and $B_{2k+1} = \{p\in[2n] \, | \, a_{2k} < p\}$, so we have
$H = \bigsqcup B_i$. We see that all the points of $H$ are fixed by $\sigma$, so for any $h\in H$, $A'(h) = A(h)$ which is not in $\{a_i\}_{i=1}^{2k}$ so the structure of this
partition is preserved in $A'$. We now group the sets $B_i$ by isolating the $k$ points $a^-_m$ with $A'(a^-_m) < 0$, constructing the $k+1$ sets $V_m$ such that $V_m$ is the union of
all the sets $B_i$ such that for each $b\in B_i$ $a^-_{m} < b < a^-_{m+1}$ where in the case $m = 0$ or $k$ the left or right bound is ignored, respectively. 

Each set $V_m$ determines
the half-edges to be placed around a vertex of $\Psi(A)$: $V_0$ gives the half-edges to be placed around the root vertex, and under the same order-isomorphic correspondence
between the points $-k, \ldots, -1$ and the non-root vertices of $\Psi(A)$, $V_m$ gives the half-edges to be placed around the vertex corresponding to $A'(a^-_{m})$ for $m>0$.

\begin{figure}[h]
 \centering
\includegraphics[scale=0.9]{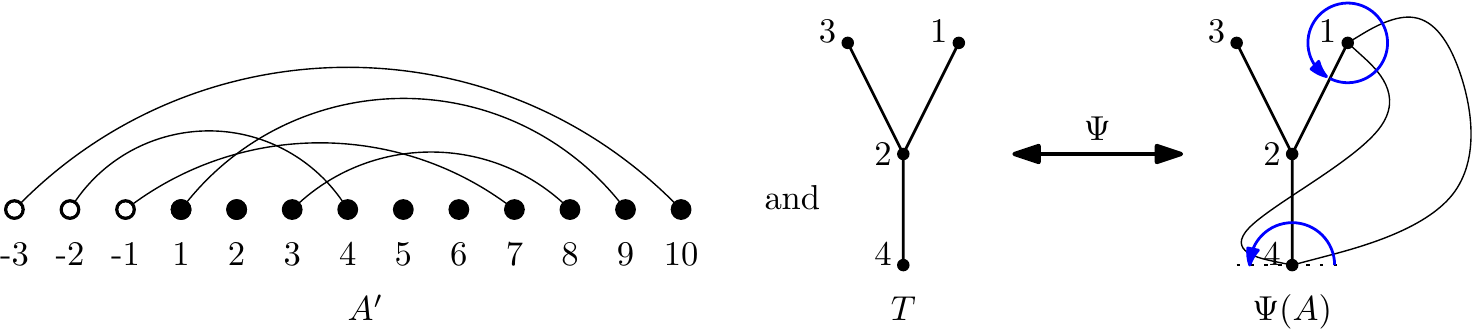}
\caption{The points 1, 3, 8, and 9 correspond to half-edges in $\Psi(A)$, which are connected the same way they are paired in 
$A'$, and in the same order around the vertices of $\Psi(A)$ (indicated by the blue arrows).}
\label{fig:bijectionExamplePart4}
\end{figure}

We see that each $V_m$ is the disjoint union of exactly as many sets $B_i$ as there are ``slots'' around the corresponding vertex (areas between outgoing edges of $T$),
with the counterclockwise
order of these slots (starting from the edge leading to the parent in $T$) corresponds to the natural order of the $B_i$. We place
$\#B_i$ half-edges into the slot corresponding to $B_i$. The counterclockwise order of these is similarly given by the natural order of the points in $B_i$.
So each point in $V_m$ corresponds to a precise
outgoing (currently un-paired) half edge from the corresponding vertex in $T$. The pairing of the half-edges is then given by the pairing of the points in $H$.
\begin{ex}\label{ex:bijectionExamplePart4}
 Continuing as in Example \ref{ex:bijectionExamplePart3}, given $A'$ and the spanning tree $T$, we construct $\Psi(A)$ by inserting
 and pairing half-edges in $T$. Here, $H = \{1,3,8,9\}$, with
$B_1 = \{1\}$, $B_2 = \{3\}$, and $B_3 = \{8,9\}$. $a^-_1 = 4$, $a^-_2 = 7$, and $a^-_3 = 10$, so $V_0 = B_1\sqcup B_2$, $V_1 = \emptyset$, $V_2 = B_3$, and $V_3 = \emptyset$.
So, 1 corresponds to the first half-edge (in counterclockwise order) around the root, 3 is the 2nd half-edge around the root but it is on the other side of $T$ because $A'(2) = 2$, 
and 8 and 9 are two half-edges around vertex 1. Because $A'(1) = 9$ and $A'(3) = 8$, the first half-edge around the root is paired
with the second half-edge around vertex 1, leaving the 2nd half-edge around the root to be paired with the 1st half-edge around vertex 1.
This is shown in Figure \ref{fig:bijectionExamplePart4}. \end{ex}

\subsection{Computing $\Psi^{-1}$}
To compute the inverse $\Psi^{-1}(M)$ for some tree-rooted map $M$ with root vertex $j+1$ and spanning tree $T$, first 
produce a Dyck path $P' = s_1's_2'\ldots s_{2k}' = \beta^{-1}(T)$ from the spanning
tree of $M$ using the inverse of the bijection $\beta$ as in Lemma~\ref{lem:pathTree}. Each down step $s$ in $P'$ 
corresponds to some vertex $v\in [k+1]\backslash\{j+1\}$ of $M$ as in Remark \ref{rem:correspondence}. The label of the down 
step $s$ is then given by element of $\{-k,-k+1,\ldots,-1\}$ corresponding to $v$ under the order-isomorphic correspondence 
between $\{-k,-k+1,\ldots,-1\}$ and $[k+1]\backslash\{j+1\}$. 

Now define an ordered sequence of points $a_i$ (the labels 
$a_i$ are currently just labels but will later take on specific integer values), setting $A(a_i) = a_i$ if $s_i'$ is an up 
step, and $A(a_i) = p$ if $s_i'$ is a down step with label $p$. We denote by $n_k$ the subsequence of indices with 
$A(a_{n_k}) < 0$ (which were down steps in $P'$). As described previously, each of the $a_{n_k}$ 
corresponds to a unique non-root vertex $v$ of $M$, and
each point $a_i$ with $n_k < i < n_{k+1}$ is a fixed point and corresponds to an outgoing half-edge from $v$ which is 
part of $T$, with the left-to-right order of these fixed points corresponding to the counterclockwise order of the outgoing 
half-edges from $v$. So each pair $(a_i,a_{i+1})$ corresponds to a unique slot between two edges of $T$. 

So, for each outgoing 
half-edge from $v$ which is \textit{not} part of $T$, we add a point between $a_i$ and $a_{i+1}$ where $(a_i,a_{i+1})$ corresponds 
to the slot containing the half-edge. If $p_1$ and $p_2$ are two points added in this manner, and the half-edges to which they 
correspond are the two halves of the same edge, then we set $A(p_1) = p_2$ (and 
therefore also $A(p_2) = p_1$, since $A$ is an involution). This brings the total number of positive points to $2n$. Now label the 
positive points in order with the integers in $[2n]$ (this also gives the definite values for the labels $a_i$). Lastly, find 
the path $P = s_1s_2\ldots s_{2k}$ with $j$ flaws which corresponds to $P'$ under the correspondence described in Lemma~\ref{lem:ChungFeller}. 
Define $\sigma$ as before, setting $\sigma(a_i) = a_j$ if and only if $s_i = s_j'$. Then $\Psi^{-1}(M)$ is given by 
$\sigma^{-1}A'\sigma$. So $\Psi$ is bijective. \qed

This concludes the proof of Theorem \ref{theor:Bijection}.

\section{Orthogonal Polynomials and $q$-analogues}\label{sec:orthoPolys}

With this combinatorial connection established, we study rook placements and related $q$-analogues to work towards
a $q$-analogue of one side of \eqref{eq:HZForm}, beginning with a few more definitions: A \textit{Motzkin path} $\omega$ is a sequence of lattice points
$(p_0 = (x_0,y_0), p_1 = (x_1,y_1), \ldots, p_n = (x_n,y_n))$ such that for each $i$, $(x_{i+1},y_{i+1}) = (x_i + 1, y_i - 1)$ or $(x_i + 1, y_i)$
or $(x_i + 1, y_i + 1)$. The pair $(p_{i-1}, p_i)$ is called the $i^{\textrm{th}}$ step of $\omega$. The $i^{\textrm{th}}$ step of $\omega$ is said to be southeast,
east, or northeast if $p_{i-1}$ and $p_i$ satisfy the first, second or third of the possible relations, respectively. The $i^{\textrm{th}}$ step is said
to start on level $y_{i-1}$ and end on level $y_i$.

In \cite{GV}, Viennot elucidated the following important result (due to Viennot and Flajolet): 
\begin{lem}[\cite{GV}]\label{lem:orthoPolyMoments} For polynomials defined by $P_{n+1} = (x-b_n)P_n - \lambda_nP_{n-1}$ with $P_0 = 1$ and $P_1 = x$, which are
orthogonal with respect to a weight function $w(x)$,
\[
 \int_a^bx^nP_iP_jw(x)dx = \lambda_1\lambda_2\ldots\lambda_j\sum_{\omega} V(\omega),
\]
where the sum on the right is over Motzkin paths $\omega$ of length $n$ starting on level $i$ and ending on level $j$, and
$V(\omega)$ is defined as the product of the weights of each step in the path, where the weight of a step starting on level 
$k$ is 1 if it is northeast, $b_k$ if it is east, and $\lambda_k$ if it is southeast.
\end{lem}
\begin{rem}\label{rem:qLambdas}
 In principle, one could define the $\lambda_n$ as polynomials in some other parameter $q$; this will be the basis for the $q$-Hermite polynomials.
\end{rem}

\begin{figure}[h]
\centering
\includegraphics[scale=0.9]{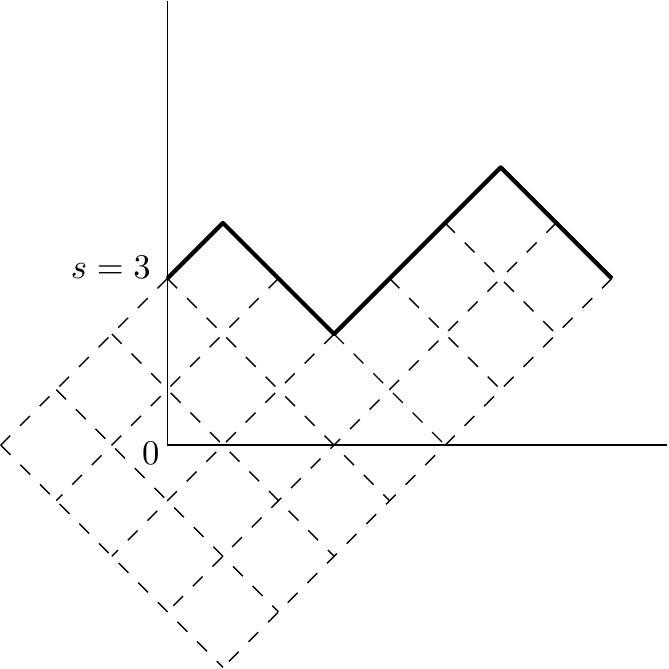}
\caption{The correspondence between weighted paths for Hermite polynomials and boards. In this case we have $s = 3$ and $\mu = (3+1,3+1,3+4,3+4)$. $s$ is always given by the starting and ending levels of the path and $2n$ is the length of the path (ensuring $n$ rows in the board). The starting level of the southeast step at the end of row $i$ is $\mu_i-i+1$ because there were $i-1$ previous southeast steps. }
\label{fig:weightedPathEx}
\end{figure}

The Hermite polynomials $H_n$ are a well-studied family of polynomials orthogonal with respect to the Gaussian measure $w(x)\textrm{d}x = \frac{1}{\sqrt{2\pi}}e^{-\frac{x^2}{2}}\textrm{d}x$ (on the real line),
defined by the recurrence $H_{n+1} = xH_n -nH_{n-1}$ with $H_0 = 1$ and $H_1 = x$. With the above definitions, $b_n = 0$ and $\lambda_n = n$ for all $n \geq 0$. Moments
of these polynomials of the form $\frac{1}{\sqrt{2\pi}}\int_{-\infty}^{\infty}x^{2n}H_s^2e^{-\frac{x^2}{2}}\textrm{d}x$ are closely related to rook placements
on boards in $Y(n,s)$: by Lemma~\ref{lem:orthoPolyMoments}, this integral is a sum over weighted Motzkin paths of length $2n$ beginning and ending at level $s$. However, because
$b_n = 0$ we can take the sum to be over only paths with no east segments. Paths with no east segments of length $2n$ beginning and ending at level $s$ are in
obvious correspondence with boards in $Y(n,s)$ (see Figure~\ref{fig:weightedPathEx}). 

Because $\lambda_n = n$, the weight $V(\omega)$ of such a path $\omega$ exactly counts the number
of placements of $n$ non-attacking rooks on the corresponding board (see (\ref{eq:rpEnum})). So we have that
\begin{equation}\label{eq:hermiteMomentRookPlacements}
 \frac{1}{\sqrt{2\pi}}\int_{-\infty}^{\infty}x^{2n}H_s^2e^{-\frac{x^2}{2}}\textrm{d}x = s!\sum_{\omega}V(\omega) = s!(\#RC(n,s)).
\end{equation}

We will also need some preliminary definitions to discuss $q$-analogues. Many $q$-analogues are built from simpler $q$-analogues, the simplest being the $q$-analogue of an integer $n$, defined as $[n]_q = 1 + q + q^2 + \ldots + q^{n-1} = \frac{1-q^n}{1-q}$.
The notation $[\cdots]_q$ is sometimes referred to as the $q$-bracket (not to be confused with the notation $[n] = \{1, \ldots, n\}$), 
and generally denotes a $q$-analogue. For example, from this simple $q$-analogue one can then define
a $q$-analogue of factorial, the $q$-factorial $[n]_q! = \prod_{i=1}^n[i]_q = [n]_q[n-1]_q\cdots[2]_q[1]_q$. One can then further define the $q$-binomial coefficient
\[
\qbinom{n}{k}_q = \frac{[n]_q!}{[k]_q![n-k]_q!} = \frac{[n]_q[n-1]_q\cdots[n-k+1]_q}{[k]_q!}
\]
and the $q$-double factorial
\[
[2n-1]_q!! = [2n-1]_q[2n-3]_q\cdots [3]_q [1]_q. 
\]

We briefly recall how the above $q$-analogues of $\binom{n}{k}$
and $(2n-1)!!$ can be seen as a generating series with a certain
statistic.

\begin{enumerate}
\item[1.] The following can be found in \cite[page 65]{EC1}:
\begin{equation} \label{qnchoosek}
\qbinom{n}{k}_q = \sum_{\alpha} q^{\sz(\alpha)},
\end{equation}
Where the sum is over all choices of a $k$-element subset $\alpha \subset [n]$, and the statistic
$\sz$ is given by 
\[
 \sz(\alpha) := \sum_{i \in \alpha} \#\{j \in [n] \mid j \not\in \alpha; \, j < i\}.
\]
Note that it follows immediately from this identity that the $q$-binomial coefficients are, in fact, polynomials in $q$ with 
nonnegative integer coefficients.

\item[2.] The following can be found in \cite[(5,4)]{SiSt} and \cite{RDEP}:
\begin{equation} \label{qcrossingnestings}
[2n-1]_q!!  = \sum_{m\in \mathcal{M}_{2n}} q^{\cn(m)},
\end{equation}
where $\mathcal{M}_{2n}$ is the set of perfect matchings on $[2n]$ and the statistic $\cn(m)$ is given by 
$\#\{{\rm crossings}\} + 2(\#\{{\rm nestings}\})$ where a crossing is an instance of arcs $(i, j)$ and $(k, l)$ with
$i < k < j < l$ and a nesting is an instance of arcs $(i, j)$ and $(k, l)$ with $i < k < l < j$.

\item[3.] Note that we could define the statistics $\sz$ and $\cn$ on perfect matchings $m$ of a $2k$-element subset 
$\alpha \subset [n]$ (strictly speaking, we would have a pair $(\alpha, m)$ of a subset and a matching, and we could then
take the matching to be of $[2k]$) by having $\sz(m) = \sz(\alpha)$ and $\cn(m)$ defined as before (the points outside 
$\alpha$ contribute nothing to the number of crossings or nestings as they are defined above, as those only depend on arcs 
of the matching), and with this we would have
\begin{equation}\label{combinedStats}
 \sum_{(\alpha, m)} q^{\sz(m) + \cn(m)} = \qbinom{n}{2k}_q[2k - 1]_q!!,
\end{equation}
where the sum is over all choices of a $2k$-element subset of $[n]$ and a matching of that subset (all the pairs $(\alpha, m)$
as described above). This identity follows from the previous two, as well as the fact that the matching and the choice of 
a subset can be made completely independently.

\end{enumerate} 

\begin{figure}[h]
\centering
\includegraphics{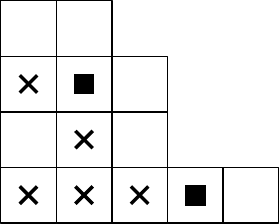}
\caption{Counting inversions. For each rook, the square where it is, and the squares directly to the left of it, and directly below it are crossed out. The remaining squares are inversions. In this example, inv $=6$.}
\label{fig:inversionsExample}
\end{figure}

In \cite{GR}, Garsia and Remmel introduced a $q$-analogue of the number of rook placements on a Young diagram-shaped board $\mu$ with $n$ rows:
\begin{equation}\label{eq:qRookCounting}
 \prod_{i=1}^n[\mu_i-i+1]_q = \sum_{\substack{n\textrm{-rook placements}\\ C \,\textrm{on } \mu}}q^{\textrm{inv}(C)},
\end{equation}
where the statistic $\inv$ (the inversions) is obtained from a board
and a rook placement on it by the following procedure: for each rook,
cross out the square where the rook is, each square below it in its column
and each square to the left of it in its row. After this has been done
for all rooks, the number of squares on $\mu$ that have not been crossed out is the number of inversions $\inv(C)$. For an example,
see Figure~\ref{fig:inversionsExample}. This identity is proven by a simple induction on the rows of the board.

To express the $q$-analogues in Conjecture~\ref{conj:qIdent} as moments of polynomials, we will also need $q$-analogues of
the Hermite polynomials, a $q$-analogue of the Gaussian measure (in particular, one by which our $q$-Hermite polynomials
are orthogonal), and a $q$-analogue of integration.

Our $q$-Hermite polynomials $H_n(x,q)$ (the symbol $H_n$ will refer to these $q$-analogues from here on) will be defined by 
the recurrence
\[
 H_{n+1}(x,q) = xH_n(x,q) - q^{n-1}[n]_qH_{n-1}(x,q), \quad H_0 = 1 \textrm{ and } H_1 = x.
\]

We note that these are closely related to the Discrete $q$-Hermite I polynomials $h_n(x,q)$ \cite{KS} defined by the recurrence
\[
 h_{n+1}(x,q) = xh_n(x,q) - q^{n-1}(1-q^n)h_{n-1}(x,q) \quad \textrm{with} \quad h_0 = 1 \quad \textrm{and} \quad  h_1 = x.
\]
The relationship is the following: 
\begin{equation}\label{eq:qhermRel}
 H_n(x,q) = a^{-n}h_n(ax,q), \quad {\rm where} \quad a = \sqrt{1-q}, 
\end{equation}
 which follows simply by checking the recurrence.
 

In fact, we have an explicit formula for $H_n$, which we will later use to set up
a ``$q$-inclusion-exclusion'' formula which could potentially be used to prove Conjecture \ref{conj:qIdent} combinatorially:
\begin{lem}\label{lem:qHermiteExplicitSum}
 We have the following formula for $H_n(x,q)$:
\[
 H_n(x,q) = \sum_{k=0}^{\lfloor n/2 \rfloor} (-1)^kq^{k(k-1)}\qbinom{n}{2k}_q[2k-1]_q!!x^{n-2k}.
\]
\end{lem}
Lemma \ref{lem:qHermiteExplicitSum} will be proved in the Appendix.

A $q$-analogue of integration (also called the Jackson integral) and of
the Gaussian measure were used in \cite{RDEP}, and we will use these with our $q$-Hermite polynomials. Many facts about Jackson
integrals and $q$-differentiation can be found in \cite{KS}. We will use the notation $\int_a^b f(x)d_qx$ for the
Jackson integral.


In general, two $q$-analogues of exponential functions are defined as
\[
 e_q^x = \sum_{n=0}^{\infty}\frac{x^n}{[n]_q!} \quad \textrm{and} \quad E_q^x = \sum_{n=0}^{\infty}q^{\binom{n}{2}}\frac{x^n}{[n]_q!}.
\]
These satisfy the following relation \cite{RDEP}: $e_q^xE_q^{-x} = 1$.
With this, the Gaussian $q$-distribution is defined as
\begin{equation}\label{eq:gaussianqDistribDef}
 w_G(x) := E_{q^2}^{\frac{-q^2x^2}{[2]_q}} = \sum_{n=0}^{\infty} \frac{q^{n(n+1)}(q-1)^n}{(1-q^2)_{q^2}^n}x^{2n},
\end{equation}
where
\[
 (a+b)_q^n = \prod_{i=0}^{n-1}(a+q^ib).
\]

The moments of $w_G(x)$ are given \cite{RDEP} by
\[
 \frac{1}{c(q)}\int_{-\nu}^{\nu}x^nw_G(x)d_qx = [n-1]_q!!,
\]
where $\nu = \frac{1}{\sqrt{1-q}}$, $[n-1]_q!! = 0$ if $n$ is odd and $[n-1]_q[n-3]_q\cdots[3]_q[1]_q$ if $n$ is even, and
\begin{equation}\label{eq:cq}
 c(q) = \int_{-\nu}^{\nu}w_G(x)d_qx = 2\sqrt{1-q} \left(\sum_{i=0}^{\infty}\frac{(-1)^iq^{2\binom{i+1}{2}}}{(1-q^{2i+1})(1-q^2)_{q^2}^i}\right).
\end{equation}

In \cite[page 118]{KS}, the orthogonality relation for
the $h_n$ is given:
\begin{equation}\label{eq:littlehOrtho}
 \int_{-1}^{1}h_m(x,q)h_n(x,q)w(x,q)d_qx = (1-q)\left(\prod_{i=1}^{n}(1-q^i)\right)\tilde c(q)q^{\binom{n}{2}}\delta_{mn},
\end{equation}
where
\[
 w(x,q) = \prod_{i=1}^{\infty}(1-q^ix)\prod_{j=1}^{\infty}(1+q^jx) \quad \textrm{and} \quad \tilde c(q) = \prod_{i=1}^{\infty}\frac{1}{(1+q^i)(1-q^{2i})}.
\]
We can now prove the orthogonality of our $q$-Hermite polynomials $H_n$ with the Gaussian $q$-distribution $w_G$, which we can
then use to apply Lemma \ref{lem:orthoPolyMoments}:
\begin{lem}\label{lem:qHermRelate}
 We have the orthogonality relation
\begin{equation}\label{eq:bigHOrtho}
\int_{-\nu}^{\nu}H_m(x,q)H_n(x,q)w_G(x)d_qx = c(q)q^{\binom{n}{2}}[n]_q!\delta_{mn},
\end{equation}
where $\nu = 1/\sqrt{1 - q}$, $c(q)$ is as in \eqref{eq:cq}, and $\delta_{mn}$ is the Kronecker delta.
\end{lem}

\noindent\textit{Proof.} Making the variable transformation $x\rightarrow ax$ in (\ref{eq:littlehOrtho}) and multiplying both sides by $a^{-(m+n)}$, (\ref{eq:littlehOrtho}) becomes
\[
 a\int_{\frac{-1}{a}}^{\frac{1}{a}}a^{-m}h_m(ax,q)a^{-n}h_n(ax,q)w(ax)d_qx = (1-q)a^{-2n}\left(\prod_{i=1}^{n}(1-q^i)\right)\tilde c(q) q^{\binom{n}{2}}\delta_{mn},
\]
where we have replaced $a^{-(m+n)}$ on the right with $a^{-2n}$ because of the presence of $\delta_{mn}$.
Using \eqref{eq:qhermRel}, the fact that $a^{-2}(1-q^i) = [i]_q$, and $a = \nu^{-1}$, this is
\[
 \int_{-\nu}^{\nu}H_m(x,q)H_n(x,q)w(ax)d_qx = \left(\sqrt{1-q}\right)[n]_q!\tilde c(q) q^{\binom{n}{2}}\delta_{mn}.
\]
It now remains to show that $w(ax) = w_G(x)$ and $\left(\sqrt{1-q}\right)\tilde c(q) = c(q)$. Given $w(ax) = w_G(x)$, the latter follows easily from the fact that we must have
$\tilde c(q) = \int_{-1}^{1}w(x)d_q(x)$ since by Lemma \ref{lem:orthoPolyMoments} (which applies because of the orthogonality), the norms
of the $h_n$ are given by $\prod_{i=1}^{n}\lambda_i = \prod_{i=1}^{n}q^{i-1}(1-q^i) = q^{\binom{n}{2}}\prod_{i=1}^{n}(1-q^i)$, from which it
follows that $c(q) = a\tilde c(q) = \left(\sqrt{1-q}\right)\tilde c(q)$ by a simple variable change. To show $w(ax) = w_G(x)$,
we note that
\[
 \frac{1}{w(x,q)} = \prod_{i=1}^{\infty}\frac{1}{(1-q^ix)}\prod_{j=1}^{\infty}\frac{1}{(1+q^jx)} = \prod_{i=0}^{\infty}\frac{1}{1-q^{2i}q^2x^2}
\]
and use the $q$-binomial theorem (\cite[page 72]{EC1}) to obtain
\[
 \frac{1}{w(x,q)} = \sum_{i=0}^{\infty}\frac{\left(q^2x^2\right)^i}{[i]_{q^2}!(1-q^2)^i} = \sum_{i=0}^{\infty}\frac{\left(\frac{q^2x^2}{(1-q^2)}\right)^i}{[i]_{q^2}!} = e_{q^2}^{\frac{q^2x^2}{1-q^2}}, 
\]
meaning $w(x,q) = E_{q^2}^{-q^2x^2/(1-q^2)}$. Making the substitution $x\rightarrow ax$ gives
\[
 w(ax,q) = E_{q^2}^{-\frac{q^2(1-q)x^2}{1-q^2}} = E_{q^2}^{-\frac{q^2x^2}{[2]_q}} = w_G(x),
\]
and the lemma is proved. \qed

We can now use the orthogonality (\ref{eq:bigHOrtho}) to prove the following theorem:
\begin{theor}\label{theor:qHermiteMoments}
The moments of the $q$-Hermite polynomials against the $q$-Gaussian distribution are given by
\begin{equation}\label{eq:qHermiteMoments}
\frac{1}{q^{\binom{s}{2}}[s]_q!c(q)}\int_{-\nu}^{\nu}x^{2n}H_s^2(x,q)w_G(x)d_qx = \sum_{\mu \in Y(n,s)}\prod_{i=1}^nq^{\mu_i - i}[\mu_i - i + 1]_q,
\end{equation}
where
\[
 \nu = \frac{1}{\sqrt{1-q}} \quad \textrm{and} \quad c(q) = 2(1-q)^{\frac{1}{2}}\sum_{m=0}^{\infty}\frac{(-1)^mq^{m(m+1)}}{(1-q^{2m+1})(1-q^2)_{q^2}^m}
\]
\end{theor}
\noindent\textit{Proof.} Because of the orthogonality (\ref{eq:bigHOrtho}), we can apply Lemma \ref{lem:orthoPolyMoments}:
\[
 \frac{1}{c(q)}\int_{-\nu}^{\nu}x^{2n}H_s^2(x,q)w_G(x)d_qx = \lambda_1\lambda_2\cdots\lambda_s\sum_{\omega}V(\omega) = q^{\binom{s}{2}}[s]_q!\sum_{\mu\in Y(n,s)}\prod_{i=1}^{n}q^{\mu_i-i}[\mu_i-i+1]_q,
\]
where the product $\lambda_1\lambda_2\cdots\lambda_s$ became $q^{\binom{s}{2}}[s]_q!$ because $\lambda_i = q^{i-1}[i]_q$ and the sum over Motzkin paths of length $2n$ was replaced
with the sum over Young shapes $\mu \in Y(n,s)$ by the same method that Kerov used to do it in the $q=1$ case. \qed

\section{Towards a $q$-analogue of one side of the Harer-Zagier Formula}\label{sec:qAnalogues}

As described in Section \ref{sec:bkgrnd}, the total number of rook placements on boards with a shape described by some $\mu$ in $Y(n,s)$
is given by $\#RC(n,s) = \sum_{\mu \in Y(n,s)}\prod_{i=1}^n (\mu_i-i+1)$. Kerov's bijection between rook placements and involutions
provided an alternative way of counting rook placements which gives the expression on the right of (\ref{eq:HZForm}). So we
have the identity

\begin{equation}\label{eq:Ident}
  \sum_{\mu\in Y(n,s)}\prod_{i=1}^n(\mu_i-i+1) = \sum_{k\geq 0}(2n-1)!!\binom{s}{k}\binom{n}{k}2^k.
\end{equation}
As was discussed in Section~\ref{sec:orthoPolys}, there is a known $q$-analogue of the left side of this equation which is obtained by summing $q^{\textrm{inv}(C)}$ over elements of $RC(n,s)$, where the
statistic inv counts the number of blank squares left on $\mu$ after crossing out any square directly to the left
or above a rook (see Figure \ref{fig:inversionsExample}). This gives
\[
 \sum_{\substack{\textrm{rook placements}\\ C \,\textrm{on } \mu}}q^{\textrm{inv}(C)} = \prod_{i=1}^n[\mu_i-i+1]_q = [V(\mu)]_q.
\]
For a $q$-analogue of (\ref{eq:Ident}), we have the following conjecture:
\begin{conj}\label{conj:qIdent}
The following identity holds:
\begin{align}
\sum_{\mu \in Y(n,s)} \prod_{i=1}^{n}q^{\mu_i-i}[\mu_i-i+1]_q &= \sum_{\mu \in Y(n,s)}\sum_{\substack{\textrm{rook placements}\\ C \,\textrm{on } \mu}}q^{\textrm{inv}(C) + |\mu| - \binom{n+1}{2}} \\
&=\sum_{k\geq 0} q^{n(s-k) + \binom{k}{2}}[2n-1]_q!!\qbinom{s}{k}_q\qbinom{n}{k}_q\prod_{i=1}^k(1 + q^{n+i}), \label{eq:qIdent}
\end{align}
where $Y(n,s)$ is the set of all Young diagrams with $n$ rows of lengths between $s$ and $n+s$.
\end{conj}

The rest of this section is devoted to proving some special cases of
Conjecture~\ref{conj:qIdent} (Section~\ref{subsec:special}), and giving a
recurrence that the left-hand-side of \eqref{eq:qIdent} satisfies (Section~\ref{sec:recurr}).

\begin{rem}
 A full proof of Conjecture~\ref{conj:qIdent} was given by Stanton \cite{DSpf}, which used techniques of hypergeometric series which can be found in \cite{GRhs}. 
 However, no combinatorial proof has yet been obtained.
\end{rem}

\subsection{Some special cases of Conjecture~\ref{conj:qIdent}: $s = 0$ and $k = n$} \label{subsec:special}
We will verify conjecture \ref{conj:qIdent} for two special cases: $s=0$ and $k=n$ (that is, if we restrict the sum to only rook placements with all $n$ rooks in the first $s$ columns).
\begin{prop}\label{prop:sEquals0}
Conjecture~\ref{conj:qIdent} holds in the case $s=0$, where it becomes
\begin{equation}\label{eq:caseS=0}
 \sum_{\mu\in n^n}\sum_{\substack{\textrm{rook placements}\\ C\textrm{ on } \mu}} q^{\textrm{inv}(C) + |\mu| - \binom{n+1}{2}} = [2n-1]_q!!.
\end{equation}
\end{prop}
\noindent\textit{Proof.} By equation (\ref{eq:qRookCounting}), 
\[
 \sum_{\mu\in n^n}\sum_{\substack{\textrm{rook placements}\\ C\textrm{ on } \mu}} q^{\textrm{inv}(C) + |\mu| - \binom{n+1}{2}} = \sum_{\mu\in n^n} \prod_{i=1}^n q^{\mu_i - i}[\mu_i-i+1]_q,
\]
and by Theorem~\ref{theor:qHermiteMoments},
\[
 \sum_{\mu\in n^n} \prod_{i=1}^n q^{\mu_i - i}[\mu_i-i+1]_q = \frac{1}{c(q)}\int_{-\nu}^{\nu}x^{2n}w_G(x)d_qx = [2n-1]_q!!,
\]
where the last step uses the fact that the $2n^{\textrm{th}}$ moment of $w_G(x)$ is $[2n-1]_q!!$ \cite{RDEP}. \qed

\begin{rem}
\begin{enumerate}
\item[(i)] If instead of the statistic $\textrm{inv}(C)+|\mu|-\binom{n+1}{2}$,
  one just considers $\textrm{inv}(C)$ in the left-hand-side
  of~\eqref{eq:caseS=0}, one obtains another interesting but less
  compact $q$-analogue of $(2n-1)!!$ known as the {\em
    Touchard-Riordan formula} \cite{T}.
\item[(ii)] In \cite{RDEP}, $[2n-1]_q!!$ is expressed as a generating series over
the same objects, placements $C$ of $n$ non-attacking rooks on a
Young shape $\mu$ inside $n^n$, but with respect to a different statistic
than $\textrm{inv}(C)+|\mu|-\binom{n+1}{2}$. Such statistic is in terms of the
{\em crossings} and {\em nestings} of the arcs of the involution corresponding to
the pair $(\mu,C)$ under Kerov's bijection $\kappa$
\cite[(5,4)]{SiSt} (see \eqref{qcrossingnestings}). For general $s$, we were unable to find a corresponding equidistributed statistic for the pairs $(\mu,C)$
in $Y(n,s) \times RC_n(n,s)$. 
\end{enumerate}
\end{rem}

\begin{prop}\label{prop:CatalanCase}
 Conjecture~\ref{conj:qIdent} holds in the case $k=n$, that is,
\begin{equation}\label{eq:CatalanCase}
 \sum_{\mu \in Y(n,s)}\sum_{\substack{\textrm{rook placements}\\ C \in RC_n(n,s)}} q^{\textrm{inv}(C) + |\mu| - \binom{n+1}{2}} = q^{n(s-n) + \binom{n}{2}}[2n-1]_q!!\qbinom{s}{n}_q\prod_{i=1}^n(1+q^{n+i}).
\end{equation}
\end{prop}
\noindent\textit{Proof.} Because inversions are squares with no rook directly to the right of or above them, every square in the last $n$ columns (the $n\times n$ Young shape)
will always be an inversion. Furthermore, the inversions in the first $s$ columns are independent of the shape of the board, so the sum factors:
\begin{align}
 \sum_{\mu \in Y(n,s)}\sum_{\substack{\textrm{rook placements}\\ C \in RC_n(n,s)}} q^{\textrm{inv}(C) + |\mu| - \binom{n+1}{2}} &= \left(\sum_{\lambda \in n^n}q^{ns + 2|\lambda| - \binom{n+1}{2}}\right)\left(\sum_{\substack{\textrm{rook placements} \\ C \subset [n]\times [s]}} q^\textrm{inv(C)}\right) \\
&= q^{ns - \binom{n+1}{2}}\left(\sum_{\lambda \in n^n}q^{2|\lambda|}\right)\left(\sum_{\substack{\textrm{rook placements} \\ C \subset [n]\times [s]}} q^\textrm{inv(C)}\right).
\end{align}
By \cite[page 65]{EC1}, the left sum is given by $\qbinom{2n}{n}_{q^2}$ and by \ref{eq:qRookCounting}, the right sum is given by $\prod_{i=1}^{n}[s-i+1]_q = \qbinom{s}{n}[n]_q!$. Furthermore, $\binom{n+1}{2} = n^2 - \binom{n}{2}$, so we have
\begin{align}
\sum_{\mu \in Y(n,s)} & \sum_{\substack{\textrm{rook placements}\\ C \in RC_n(n,s)}} q^{\textrm{inv}(C) + |\mu| - \binom{n+1}{2}} = q^{n(s-n)+\binom{n}{2}}\qbinom{2n}{n}_{q^2}\qbinom{s}{n}[n]_q! \\
&= q^{n(s-n) + \binom{n}{2}}\frac{[2n]_{q^2}[2n-1]_{q^2}\cdots[n+1]_{q^2}}{[n]_{q^2}!}\qbinom{s}{n}[n]_q! \\
&= q^{n(s-n) + \binom{n}{2}}\frac{[2n]_q(1+q^{2n})[2n-1]_q(1+q^{2n-1})\cdots[n+1]_q(1+q^{n+1})}{(1+q)^n[n]_{q^2}!}\qbinom{s}{n}[n]_q! \\
&= q^{n(s-n) + \binom{n}{2}}[2n-1]_q[2n-3]_q\cdots[3]_q[1]_q\qbinom{s}{n}\prod_{i=1}^{n}(1+q^{n+i}),
\end{align}
and the proposition is proved. \qed

Because we already have the equality (\ref{eq:qHermiteMoments}), proof of Conjecture~\ref{conj:qIdent} is reduced
to proving that the integral on the left of (\ref{eq:qHermiteMoments}) also evaluates to (\ref{eq:qIdent}).

\subsection{A recurrence for Conjecture~\ref{conj:qIdent}} \label{sec:recurr}

The motivation behind the following lemma is that one could, in principle, prove Conjecture~\ref{conj:qIdent} by verifying the 
recurrence given in the lemma for (\ref{eq:qIdent}) as well and show that the integral in (\ref{eq:qHermiteMoments}) and 
(\ref{eq:qIdent}) are equal
in the $s=0$ and $s=1$ cases, the latter part of which is trivially easy. However, we have not been able to 
verify the recurrence
for (\ref{eq:qIdent}).

\begin{lem}\label{lem:intRecursion}
 Define
\[
 f(n,s) := \frac{1}{c(q)}\int_{-\nu}^{\nu}x^{2n}H_s^2(x,q)w_G(x)d_qx.
\]
 With $c(q)$ as in \eqref{eq:cq}. Then
\[
 f(n,s) = f(n+1,s-1) - 2\sum_{l=0}^{s-2}(-1)^{s-l}\frac{q^{\binom{s-1}{2}}[s-1]_q!}{q^{\binom{l}{2}}[l]_q!}f(n+1,l) + q^{2(s-2)}[s-1]_q^2f(n,s-2).
\]
With $f(n,0) = [2n-1]_q!!$ and $f(n,1) = [2n+1]_q!!$.
\end{lem}
\noindent\textit{Proof.} First we verify the initial conditions: $f(n,0) = \frac{1}{c(q)}\int_{-\nu}^{\nu}x^{2n}w_G(x)d_qx = [2n-1]_q!!$ and
$f(n,1) = \frac{1}{c(q)}\int_{-\nu}^{\nu}x^{2n+2}d_qx = [2n+1]_q!!$ both follow from the computation of the moments in \cite{RDEP}.
To verify the recurrence, first define
\[
 f^{\times}(n,s) := \frac{1}{c(q)}\int_{-\nu}^{\nu}x^{2n+1}q^{s-2}[s-1]_qH_{s-1}(x,q)H_{s-2}(x,q)w_G(x)d_qx.
\]
Using the recurrence defining $H_s(x,q)$, we have
\begin{align}
 f(n,s) &= \frac{1}{c(q)}\int_{-\nu}^{\nu}x^{2n}(xH_{s-1}(x,q)-q^{s-2}[s-1]_qH_{s-2}(x,q))^2w_G(x)d_qx \\
&= f(n+1,s-1) - 2f^{\times}(n,s) + q^{2(s-2)}[s-1]_q^2f(n,s-2).
\end{align}
Again using this recurrence, we also have
\begin{align}
 f^{\times}(n,s) &= \frac{1}{c(q)}\int_{-\nu}^{\nu}x^{2n+1}q^{s-2}[s-1]_q(xH_{s-2} - q^{s-3}[s-2]_qH_{s-3})H_{s-2}w_G(x)d_qx \\
&= q^{s-2}[s-1]_q(f(n+1,s-2)-f^{\times}(n,s-1)).
\end{align}
We now show that
\begin{equation}\label{eq:fCross}
f^{\times}(n,s) = \sum_{l=0}^{s-2}(-1)^{s-l}\frac{q^{\binom{s-1}{2}}[s-1]_q!}{q^{\binom{l}{2}}[l]_q!}f(n+1,l)
\end{equation}
by showing that the expression on the right also satisfies this recurrence and agrees with $f^{\times}$ on initial values. We have that
\begin{align}
 q^{s-2} & [s-1]_q\left(f(n+1,s-2)-\sum_{l=0}^{s-3}(-1)^{s-1-l}\frac{q^{\binom{s-2}{2}}[s-2]_q!}{q^{\binom{l}{2}}[l]_q!}f(n+1,l)\right) \\
&= q^{s-2}[s-1]_qf(n+1,s-2)+\sum_{l=0}^{s-3}(-1)^{s-l}\frac{q^{\binom{s-1}{2}}[s-1]_q!}{q^{\binom{l}{2}}[l]_q!}f(n+1,l) \\
&= \sum_{l=0}^{s-2}(-1)^{s-l}\frac{q^{\binom{s-1}{2}}[s-1]_q!}{q^{\binom{l}{2}}[l]_q!}f(n+1,l).
\end{align}
Furthermore, $f^{\times}(n,0) = 0$ since $H_{-1} = 0$, and we see that the sum on the right of (\ref{eq:fCross}) is trivially 0
because $s-2 = -2 < 0$ and $l \geq 0$. So the lemma is proved. \qed

\begin{rem} 
Another potential method of proof of Conjecture~\ref{conj:qIdent} would be to compute the integral in \eqref{eq:qHermiteMoments} explicitly using the expression
in Lemma \ref{lem:qHermiteExplicitSum} and knowledge of the moments of $w_G(x)$ to obtain the following:

\begin{multline} \label{eq:qIncEx}
 \frac{1}{q^{\binom{s}{2}}[s]_q!c(q)}  \int_{-\nu}^{\nu}x^{2n}H_s^2(x,q)w_G(x)d_qx \\
 = \sum_{t = 0}^s (-1)^t\sum_{a = 0}^t q^{a^2 + (t-a)^2 -
   t}\qbinom{s}{2a}_q\qbinom{s}{2(t-a)}_q[2a-1]_q!![2(t-a)
 -1]_q!![2(n+s-t)-1]_q!!.
\end{multline}
Then using the expansions \eqref{qcrossingnestings} and \eqref{combinedStats} for the $q$-analogues
of $(2n-1)!!$ and $\binom{n}{2k}(2k-1)!!$, \eqref{eq:qIncEx} becomes
\begin{multline} \label{eq:qinclexcl}
\frac{1}{q^{\binom{s}{2}}[s]_q!c(q)}
\int_{-\nu}^{\nu}x^{2n}H_s^2(x,q)w_G(x)d_qx =\\ 
\sum_{t=0}^s (-1)^t\sum_{a=0}^t q^{a^2 + (t-a)^2 -t}\left(\sum_{(\alpha,m_1)} q^{\cn(m_1) + \sz(\alpha)}\right)\left(\sum_{(\beta,m_2)} q^{\cn(m_2) + \sz(\beta)}\right)\left(\sum_{m_3} q^{\cn(m_3)}\right),
\end{multline}
where $m_1$ is a perfect matching on a set $\alpha\subseteq [s]$ of
size $2a$; $m_2$ is a perfect matching on a set $\beta \subseteq
\{s+1,\ldots, 2s\}$ of size $2(t-a)$; and $\gamma$  a perfect matching on a subset of size $2a$ of
$[n]$, $\beta$ is a perfect matching on a subset of size $2(t-a)$; and
$m_3$ is a perfect matching on the remaining set $[2n+2s] \backslash
(\alpha \cup \beta)$. The RHS of \eqref{eq:qinclexcl} has the flavor
of an {\em inclusion-exclusion} expression, and suggests a possible approach to a 
combinatorial proof of Conjecture~\ref{conj:qIdent} (no such proof has yet been found).
\end{rem}



\section{Concluding remarks}\label{sec:conclusion}

In Section~\ref{sec:Bijection}, we successfully constructed a bijection between tree-rooted maps and
the rook placements. However, part of the motivation behind this effort was the intent
of obtaining a bijection which preserved interesting properties of
either type of object (tree-rooted maps or rook placements). One
relevant statistic for tree-rooted maps is the degree of the vertices
(because of the {\em symmetry property} \cite[Thm. 1.3.]{BM}). We
could not explicitly trace the degree through the bijection. On the
other side,
the statistic $\inv$ counting inversions is very important in the analysis of rook placements, so if inversions could be traced through the bijection,
we could find an analogous statistic on tree-rooted maps. Moreover, tracing
inversions could also be useful since this statistic is one of the
building blocks for Conjecture~\ref{conj:qIdent}.


Though we have been unable to prove Conjecture~\ref{conj:qIdent}, we
have verified it empirically up to $n=10$ and $s=5$. If it is proved,
then there are a number of questions that arise. 

\begin{questions}
\begin{itemize}
\item[(i)] First, the expression on the right of
(\ref{eq:qIdent}) is not directly a $q$-analogue of $C(n,N)$, but
rather one of $\#RC(n,s)$. To obtain a $q$-analogue of $C(n,N)$,
one would have to take some kind of sum over the parameter $s$ from $0$ to $N-1$. However, simply the expression in (\ref{eq:qIdent})
does not give very intelligible results. One must take the sum
\[
 [C(n,N)]_q := \sum_{s=0}^{N-1}q^{-(n-1)s}[\#RC(n,s)]_q = [2n-1]_q!!\sum_{k\geq 0}q^{\binom{k}{2}-(n-1)k}\qbinom{N}{k+1}_q\qbinom{n}{k}_q\prod_{i=1}^{k}(1+q^{n+i}),
\]
where $[\#RC(n,s)]_q$ denotes the expression in (\ref{eq:qIdent}). We
can ask many questions about this though:

\item[(ii)]Is it possible to emulate the basis change in \eqref{eq:HZForm}
($N^v \mapsto \binom{N}{k+1}$) in the reverse direction for the $q$-analogue? This would
ostensibly give a $q$-analogue of the numbers
$\varepsilon_g(n)$. 

\item[(iii)] Can one continue to ``reverse engineer'' a matrix integral in the $q$
case? We were able to go back as far as writing $[\#RC(n,s)]_q$ as a
moment of some $q$-Hermite polynomials (Theorem~\ref{theor:qHermiteMoments}), but it may
be possible to continue to go backwards through the process in Kerov's
paper to arrive at a matrix integral. We note that the work of Venkataramana~\cite[Section 6]{PV} 
may provide insight in this direction.
\end{itemize}
\end{questions}

\section{Appendix}
\noindent\textit{Proof of Lemma \ref{lem:qHermiteExplicitSum}.} We begin by considering the generating function
\[
 h(x,r) = \sum_{n=0}^{\infty}\frac{h_n(x,q)}{[n]_q!(1-q)^n}r^n,
\]
which is given in \cite[page 118]{KS}:
\[
 h(x,r) = \frac{\prod_{k=0}^{\infty}(1+rq^k)\prod_{k=0}^{\infty}(1-rq^k)}{\prod_{k=0}^{\infty}(1-rxq^k)} = \prod_{k=0}^{\infty}\left(\frac{1-r^2q^{2k}}{1-rxq^k}\right).
\]
By the $q$-binomial theorem (\cite[page 72]{EC1}),
\[
 \prod_{k=0}^{\infty}\left(\frac{1}{1-rxq^k}\right) = \sum_{k=0}^{\infty}\frac{(rx)^k)}{[k]_q!(1-q)^k}
\]
and
\[
 \prod_{k=0}^{\infty}(1-r^2q^{2k}) = \sum_{k=0}^{\infty}\frac{q^{k(k-1)}(-r^2)^k}{[k]_{q^2}!(1-q^2)^k},
\]
so
\begin{align}
 h(x,r) &= \left(\sum_{k=0}^{\infty}\frac{(rx)^k}{[k]_q!(1-q)^k}\right)\left(\sum_{k=0}^{\infty}\frac{q^{k(k-1)}(-r^2)^k}{[k]_{q^2}!(1-q^2)^k}\right) \\
&= \sum_{k=0}^{\infty}\sum_{j=0}^{\infty}\frac{q^{j(j-1)}x^k(-1)^jr^{k+2j}}{[k]_q!(1-q)^k[j]_{q^2}!(1-q^2)^j} \\
&= \sum_{n=0}^{\infty}\left(r^n\sum_{j=0}^{\lfloor n/2 \rfloor}\frac{q^{j(j-1)}x^{n-2j}(-1)^j}{[n-2j]_q!(1-q)^{n-2j}[j]_{q^2}!(1-q^2)^j}\right).
\end{align}
We now note that since $H_n(x,q) = a^{-n}h_n(ax,q)$, if we define
\[
 H(x,r) = \sum_{n=0}^{\infty}\frac{H_n(x,q)}{[n]_q!}r^n,
\]
we see that $H(x,r) = h(ax,ar)$, and since $H_n(x,q)$ is $[n]_q!$ times the coefficient of $r^n$ in $H(x,r)$, we have
\begin{align}
 H_n(x,q) &= [n]_q!a^n\sum_{k=0}^{\lfloor n/2 \rfloor}\frac{q^{k(k-1)}(ax)^{n-2k}(-1)^k}{[n-2k]_q!(1-q)^{n-2k}[k]_{q^2}(1-q^2)^{k}} \\
&= \sum_{k=0}^{\lfloor n/2 \rfloor}\frac{[n]_q!q^{k(k-1)}x^{n-2k}(-1)^k}{[n-2k]_q![k]_{q^2}!(1+q)^k} \\
&= \sum_{k=0}^{\lfloor n/2 \rfloor}\frac{[n]_q!q^{k(k-1)}x^{n-2k}(-1)^k}{[n-2k]_q![2k]_q!}[2k-1]_q!! \\
&= \sum_{k=0}^{\lfloor n/2 \rfloor}(-1)^kq^{k(k-1)}\qbinom{n}{2k}_q[2k-1]!!x^{s-2k},
\end{align}
and the lemma is proved. \qed

\bibliographystyle{plain}
\bibliography{biblio-SPUR}

\medskip

\noindent Max Wimberley, Massachusetts Institute of Technology, Cambridge, MA USA 02139, \\
maxwimberley@gmail.com

\end{document}